\documentclass{amsart}
\usepackage{amssymb,amsmath, euscript}
\input{epsf}

\newcounter{sec}

\newcounter{punct}[sec]

\def\punct{\refstepcounter{punct}{\arabic{sec}.\arabic{punct}.  }}

\newtheorem{theorem}{Theorem}[sec]
\newtheorem{proposition}[theorem]{Proposition}

\newtheorem{lemma}[theorem]{Lemma}

\newtheorem{corollary}[theorem]{Corollary}

\def\COUNTERS{\addtocounter{sec}{1}
              \setcounter{punct}{0}
          \setcounter{equation}{0}
          \setcounter{theorem}{0}
          }
          
          \def\sm{\smallskip}

\begin{document}

\newcommand{\rk}{\mathop {\mathrm {rk}}\nolimits}
\newcommand{\Aut}{\mathop {\mathrm {Aut}}\nolimits}
\newcommand{\Out}{\mathop {\mathrm {Out}}\nolimits}
\renewcommand{\Re}{\mathop {\mathrm {Re}}\nolimits}

\def\frg{\mathfrak g}
\def\frh{\mathfrak h}
\def\fin{\mathrm {fin}}

\def\ov{\overline}

\def\wt{\widetilde}

\renewcommand{\rk}{\mathop {\mathrm {rk}}\nolimits}
\renewcommand{\Aut}{\mathop {\mathrm {Aut}}\nolimits}
\renewcommand{\Re}{\mathop {\mathrm {Re}}\nolimits}
\renewcommand{\Im}{\mathop {\mathrm {Im}}\nolimits}
\newcommand{\coker}{\mathop {\mathrm {coker}}\nolimits}

\def\bfa{\mathbf a}
\def\bfb{\mathbf b}
\def\bfc{\mathbf c}
\def\bfd{\mathbf d}
\def\bfe{\mathbf e}
\def\bff{\mathbf f}
\def\bfg{\mathbf g}
\def\bfh{\mathbf h}
\def\bfi{\mathbf i}
\def\bfj{\mathbf j}
\def\bfk{\mathbf k}
\def\bfl{\mathbf l}
\def\bfm{\mathbf m}
\def\bfn{\mathbf n}
\def\bfo{\mathbf o}
\def\bfp{\mathbf p}
\def\bfq{\mathbf q}
\def\bfr{\mathbf r}
\def\bfs{\mathbf s}
\def\bft{\mathbf t}
\def\bfu{\mathbf u}
\def\bfv{\mathbf v}
\def\bfw{\mathbf w}
\def\bfx{\mathbf x}
\def\bfy{\mathbf y}
\def\bfz{\mathbf z}

\def\bfA{\mathbf A}
\def\bfB{\mathbf B}
\def\bfC{\mathbf C}
\def\bfD{\mathbf D}
\def\bfE{\mathbf E}
\def\bfF{\mathbf F}
\def\bfG{\mathbf G}
\def\bfH{\mathbf H}
\def\bfI{\mathbf I}
\def\bfJ{\mathbf J}
\def\bfK{\mathbf K}
\def\bfL{\mathbf L}
\def\bfM{\mathbf M}
\def\bfN{\mathbf N}
\def\bfO{\mathbf O}
\def\bfP{\mathbf P}
\def\bfQ{\mathbf Q}
\def\bfR{\mathbf R}
\def\bfS{\mathbf S}
\def\bfT{\mathbf T}
\def\bfU{\mathbf U}
\def\bfV{\mathbf V}
\def\bfW{\mathbf W}
\def\bfX{\mathbf X}
\def\bfY{\mathbf Y}
\def\bfZ{\mathbf Z}

\def\frD{\mathfrak D}
\def\frL{\mathfrak L}

\def\l{(\!(}
\def\r{)\!)}

\def\bfw{\mathbf w}

\def\R {{\mathbb R }}
 \def\C {{\mathbb C }}
  \def\Z{{\mathbb Z}}
  \def\H{{\mathbb H}}
  \def\O{\mathbb O}
\def\K{{\mathbb K}}
\def\N{{\mathbb N}}
\def\Q{{\mathbb Q}}
\def\A{{\mathbb A}}

\def\T{\mathbb T}
\def\P{\mathbb P}

\def\G{\mathbb G}

\def\cD{\mathcal D}
\def\cL{\mathcal L}
\def\cK{\mathcal K}

\def\bbA{\mathbb A}
\def\bbB{\mathbb B}
\def\bbD{\mathbb D}
\def\bbE{\mathbb E}
\def\bbF{\mathbb F}
\def\bbG{\mathbb G}
\def\bbI{\mathbb I}
\def\bbJ{\mathbb J}
\def\bbL{\mathbb L}
\def\bbM{\mathbb M}
\def\bbN{\mathbb N}
\def\bbO{\mathbb O}
\def\bbP{\mathbb P}
\def\bbQ{\mathbb Q}
\def\bbS{\mathbb S}
\def\bbT{\mathbb T}
\def\bbU{\mathbb U}
\def\bbV{\mathbb V}
\def\bbW{\mathbb W}
\def\bbX{\mathbb X}
\def\bbY{\mathbb Y}

\def\epsilon{\varepsilon}
\def\phi{\varphi}
\def\le{\leqslant}
\def\ge{\geqslant}

\def\B{\mathrm B}

\def\GL{\mathrm{GL}}
\def\Sp{\mathrm{Sp}}
\def\U{\mathrm{U}}
\def\dist{\mathrm{dist}}

\def\cO{\mathcal O}
\def\cN{\mathcal N}

\def\la{\langle}
\def\ra{\rangle}

\begin{center}
\bf \large
 On the group of infinite $p$-adic matrices with integer elements
 
 \medskip
 
 \sc
 Yury A. Neretin%
 \footnote{The research was carried out at the IITP RAS at the expense of the Russian 
Foundation for Sciences (project № 14-50-00150).}
\end{center}

\bigskip

{\small Let $G$ be an infinite-dimensional real classical group containing
the complete unitary group (or complete orthogonal group) as a subgroup.
Then $G$ generates a category of double cosets (train) and any unitary representation
of $G$ can be canonically extended to the train. We prove a technical lemma
about the complete group $\GL$ of infinite $p$-adic matrices with integer
coefficients, this lemma implies that the phenomenon of automatic extension
of unitary representations to  trains  is valid for infinite-dimensional $p$-adic groups.
}

\section{The statement}

\COUNTERS

{\bf\punct Notation.}
Denote by $\Q_p$ the field of $p$-adic numbers, by $\O_p$ the ring of $p$-adic integers. 
We consider 
infinite matrices $g=\{g_{ij}\}$, where $i$, $j\in\N$, over $\O_p$.
We define 3 versions of the group $\GL(\infty)$ over $\O_p$.

\sm

1) Our main object is the group $\GL(\infty,\O_p)$,
which consists of invertible matrices satisfying two conditions: 

\sm

$A^*$. for each 
$i$ we have $\lim_{j\to\infty}|g_{ij}|=0$;

\sm

$B^*$. for each $j$ we have $\lim_{i\to\infty}|g_{ij}|=0$.

\sm

2) We also consider a larger group $\ov\GL(\infty,\O_p)$
consisting of invertible matrices satisfying condition $A^*$.

\sm


3) We regard compact groups $\GL(n,\O_p)$ as subgroups of $\GL(\infty,\O_p)$ 
consisting of block $(n+\infty)$-matrices
of the form
$\begin{pmatrix}*&0\\0&1\end{pmatrix}$.

We say that an infinite matrix $g$ is {\it finitary} if $g-1$ has only finite number of nonzero matrix elements.
Denote by $\GL_\fin(\infty, \O_p)$ the group of invertible finitary infinite matrices over $\O_p$,
this group is an inductive limit
$$
\GL_\fin(\infty, \O_p)=\lim_{\longrightarrow} \GL(n,\O_p)
$$
and is equipped with a topology of an inductive limit: a function on $\GL_\fin(\infty, \O_p)$
is continuous iff its restriction to each prelimit subgroup is continuous.

\sm

{\sc Remark.} The group $\ov\GL(\infty,\O_p)$ appears in the context of \cite{Ner-hua}.
However, $\GL(\infty,\O_p)$ is a more interesting object from a point of view of \cite{Ner-p}.

\sm

{\bf\punct The result of the paper.}
Denote by $\theta_j$ the following matrix
$$
 \theta_j=
\begin{pmatrix}
 0&1_j&0\\
 1_j&0&0\\
 0&0&1_\infty
\end{pmatrix}\in \GL(\infty,\O_p),
$$
where $1_j$ denotes the unit matrix of size $j$.
The purpose of this note is to prove the following statement:

\begin{lemma}
\label{l:1}
 Consider a unitary representation $\rho$ of the group $\GL(\infty,\O_p)$
 in a Hilbert space $H$. Denote by $H^\GL\subset H$ the space of all vectors 
 fixed by all operators $\rho(g)$.
 Then the sequence
 $\rho(\theta_j)$ weakly converges to the orthogonal projection to $H^\GL$.
\end{lemma}

Since $\GL(\infty,\O_p)$ is dense in $\ov\GL(\infty,\O_p)$, we get the following corollary.

\begin{corollary}
 The same statement holds for the group $\ov\GL(\infty,\O_p)$.
\end{corollary}

{\bf \punct Variations.}
Define the orthogonal group $\mathrm{O}(\infty,\O_p)$
as the subgroup in $\GL(\infty,\O_p)$ consisting of all matrices
satisfying $g^tg=1$, where $^t$ denotes the transposing.
Denote by $J$ the $2\times 2$-matrix 
$\begin{pmatrix} 0&1\\-1&0\end{pmatrix}$ over $\O_p$. 
Denote 
\begin{equation}
I:=
\begin{pmatrix}
 J&0&\dots\\
 0&J&\dots\\
 \vdots&\vdots&\ddots 
\end{pmatrix}.
\label{eq:I}
\end{equation}
Denote by
$\mathrm{Sp}(\infty,\O_p)$ the subgroup in $\GL(\infty,\O_p)$ consisting of all matrices
satisfying $g^t I g=1$.

\sm

{\it Lemma} \ref{l:1}   ({\it with the same proof}) {\it holds for the groups}
$\mathrm{O}(\infty,\O_p)$, $\mathrm{Sp}(\infty,\O_p)$;
for $\mathrm{Sp}(\infty,\O_p)$ we must consider the sequence $\theta_{2m}\in \mathrm{Sp}(\infty,\O_p)$.

\sm 

{\bf\punct Admissibility in Olshanski's sense.} We also prove the following
technical statement.
Consider a unitary representation $\rho$ of the group $\GL_\fin(\infty,\O_p)$
in a Hilbert space $H$. Denote by $H_m$ the space of $\GL^{[m]}_\fin(\infty,\O_p)$-invariant vectors. 
We say that a representation $\rho$ is {\it admissible} (see \cite{Olsh-semigr})
if the subspace $\cup_{m=0}^\infty H_m$ is dense in $H$ .

\begin{lemma}
\label{l:adm}
 The following conditions for a representation $\rho$ of $\GL_\fin(\infty,\O_p)$ are equivalent:
 
 \sm
 
 $\bullet$ The representation $\rho$  admits a continuous extension to $\GL(\infty,\O_p)$.
 
 \sm
 
  $\bullet$  The representation $\rho$ is admissible.
\end{lemma}

{\bf\punct Structure of the paper.}
Lemma 
\ref{l:1} seems rather technical,
however it implies that $\GL(\infty,\O_p)$ is a heavy group
in the sense of \cite{Ner-book}, Chapter VIII. This implies 
numerous 'multiplicativity theorems', an example is discussed in the next section.
Lemma \ref{l:1} is proved in Section 3, Lemma \ref{l:adm} in Section 4.

\section{Introduction. An example of  multiplicativity theorems}

\COUNTERS

{\bf \punct Initial data.%
\label{ss:initial}}  
Denote by $S_\fin(\infty)$ the group of all finitely supported permutations of $\N$. Fix a ring $R$.
Let $G$ be a subgroup in $\GL_\fin(\infty,R)$,  $K$  its subgroup. Assume that $K$
contains $S_\fin(\infty)$ embedded as the group of all $0-1$ matrices.

\sm

{\sc Examples.} a) $G=K=S_\fin(\infty)$.

\sm

b) $G=\GL_\fin(\infty,\R)$, $K=\mathrm{O}_\fin(\infty)$

\sm

c) Let $R$ be the algebra of $2\times 2$ real matrices. Let $J=\begin{pmatrix}
           0&1\\-1&0
          \end{pmatrix}\in R$. We consider the group $G=\Sp_\fin(2\infty,\R)$
          consisting of matrices over $R$ preserving the skew-symmetric bilinear form
          with matrix $I$ given by (\ref{eq:I}).
          The subgroup $K=\U_\fin(\infty)$ consists of matrices whose entries have form
          $\begin{pmatrix}
            a&b\\-b&a
           \end{pmatrix}\in R$. 
    
    \sm
    
 d) $G=\GL_\fin(\infty,\Q_p)$, $K=\GL_\fin(\infty,\O_p)$.

 \sm
 
 {\sc Remark.} Denote by $G(\alpha)$ (resp. $K(\alpha)$) the subgroup in $G$ (resp. $K$)
 consisting of all $(\alpha+\infty)$-block matrices of the form
 $\begin{pmatrix} w&0\\0&1\end{pmatrix}$. These groups contain at least the finite symmetric group
 $S(\alpha)$. Then
 \begin{equation}
 G=\lim_{\longrightarrow} G(\alpha),\qquad  K=\lim_{\longrightarrow} K(\alpha).
 \label{eq:GK-limits}
 \end{equation}

{\bf \punct The multiplication of double cosets.}
We fix $n$ and consider the product $\wt G$ of $n$ copies of the group $G$,
$$\wt G= G\times G\times  \dots \times G.$$
We write elements of this product by
\begin{equation}
g=\{g^{(l)}\}:=(g^{(1)},\dots, g^{(n)}), \qquad \text{where $g_j\in G$.}
\label{eq:g()}
\end{equation}
Consider the diagonal subgroup $K\subset \wt G$, i.e., the group, whose elements are
collections
$$
(u,\dots,u),\qquad\text{where $u\in K$.}
$$

Let $\alpha=0$, $1$, $2$, \dots.
Denote by $K^\alpha$ the subgroup of $K$ consisting of all matrices
having the form $\begin{pmatrix}
                  1_\alpha&0\\0&u
                 \end{pmatrix}\in K$. Denote by
                 $$K^\alpha\setminus\wt G/K^\alpha$$
the double cosets, i.e., the space of collections (\ref{eq:g()}) defined up to the equivalence
$$
\bigl(g^{(1)},\dots,g^{(n)}\bigr)\sim \Biggl(\begin{pmatrix}
                                     1_\alpha&0\\0&u
                                    \end{pmatrix} g^{(1)} \begin{pmatrix}
                                     1_\beta&0\\0&v
                                    \end{pmatrix}, \dots, \begin{pmatrix}
                                     1_\alpha&0\\0&u
                                    \end{pmatrix} g^{(n)} \begin{pmatrix}
                                     1_\beta&0\\0&v
                                    \end{pmatrix}\Biggr),
$$
where $u$, $v\in K$.

For each $\alpha$ we define the sequence $\theta_j^{[\alpha]}$ by
$$
\begin{pmatrix}
 1_\alpha&0&0&0\\
 0&0&1_j&0\\
 0&1_j&0&0\\
 0&0&0&1_\infty
\end{pmatrix}\in K^\alpha\cap S_\fin(\infty).
$$

The following statements a)--c) can be verified in a  straightforward way
(see a formal proof in \cite{GN} for $G=S_\fin(\infty)$,
 which is valid in  a general case):

\sm

a) {\it Let $\frg_1\in K^\alpha\setminus \wt G/K^\beta$, $\frg_2\in K^\beta\setminus \wt G/K^\gamma$.
Let $g_1$, $g_2\in \wt G$ be their representatives. Then the sequence 
$$
K^\alpha g_1 \theta_j^{[\beta]} g_2 K^\gamma\in K^\alpha\setminus  \wt G/K^\gamma
$$
of double cosets
is eventually constant. Moreover the limit does not depend  on
the choice of representatives $g_1\in\frg_1$, $g_2\in\frg_2$.
}

\sm

b) {\it Thus we get a multiplication $(\frg_1,\frg_2)\mapsto\frg_1\circ\frg_2$
$$
K^\alpha\setminus \wt G/K^\beta\, \times\, K^\beta\setminus \wt G/K^\gamma \to 
K^\alpha\setminus \wt G/K^\gamma,
$$
which can be described in the following way.
We write representatives $g_1\in \frg_1$, $g_2\in \frg_2$
as block $(\alpha+\infty)\times  (\beta+\infty)$ and collections of
$(\beta+\infty)\times (\gamma+\infty)$-%
matrices  
$$
\{g_1^{(l)}\}=\left\{\begin{pmatrix}
   a_1^{(l)}&b_1^{(l)}\\c_1^{(l)}&d_1^{(l)}
  \end{pmatrix}\right\},\qquad
  g_2^{(l)}=\left\{
\begin{pmatrix}
   a_2^{(l)}&b_2^{(l)}\\c_2^{(l)}&d_2^{(l)}
  \end{pmatrix}\right\},
$$
then
a representative of $\frg_1\circ\frg_2$ is given by 
 \begin{equation}
 \{(g_1\circledcirc g_2)^{(l)}\}:=
 \left\{ \begin{pmatrix}
a_1^{(l)}&b_1^{(l)}&0\\
c_1^{(l)}&d_1^{(l)}&0\\
0&0&1_\infty
       \end{pmatrix}   
       \begin{pmatrix}
a_2^{(l)}&0&b_2^{(l)}\\
0&1_\infty&0\\
c_2^{(l)}&0&d_2^{(l)}
       \end{pmatrix}
 \right\}.
  \label{eq:uzly}      
 \end{equation}
}

The size of these matrices is
$$
\bigl(\alpha+[\infty+\infty]\bigr)\times \bigl(\gamma+[\infty+\infty]\bigr)=
\bigl(\alpha+\infty\bigr)\times \bigl(\gamma+\infty\bigr),
$$
so this collection can be regarded as a representative of  an element of the space
$K^\alpha\setminus \wt G/K^\gamma$. More precisely, we must choose  arbitrary
bijections $\sigma_1$, $\sigma_2$ between a disjoint union $\N\coprod \N$ and $\N$ to get an element of
a desired size:
$$
\left\{\begin{pmatrix}
 1_\alpha&0\\0&\sigma_1
\end{pmatrix}
 \Bigl( g_1\circledcirc g_2\Bigr)^{(l)}
 \begin{pmatrix}
 1_\alpha&0\\0&\sigma_2
\end{pmatrix}^{-1}\right\}.
$$
The double coset containing this matrix does not depend on a choice
of $\sigma_1$, $\sigma_2$.

\sm

c) {\it The product {\rm of double cosets} is associative, i.e., for 
$$
\frg_1\in 
K^\alpha\setminus \wt G/K^\beta, \qquad \frg_2\in K^\beta\setminus \wt G/K^\gamma, \qquad 
\frg_3\in K^\gamma\setminus \wt G/K^\delta
$$
we have
$$
(\frg_1\circ\frg_2)\circ\frg_3=\frg_1\circ(\frg_2\circ\frg_3).
$$
}

{\sc Remark.} The formula for  the $\circledcirc$-product 
$$
\begin{pmatrix}
 a_1&b_1\\c_1&d_1
\end{pmatrix}
\circledcirc
\begin{pmatrix}
 a_2&b_2\\c_2&d_2
\end{pmatrix}
:=
\begin{pmatrix}
 a_1 a_2& b_1& a_1 b_2\\
 c_1 a_2& d_1& c_1 b_2\\
 c_2&0&d_2
\end{pmatrix}
$$
of matrices
initially arose as a formula for
a product of operator colligations, see \cite{Brod1}, \cite{Brod}.
\hfill $\boxtimes$

 \sm
 
 {\bf \punct Multiplicativity theorems.}
Next, consider a unitary representation $\rho$ of $\wt G$
in a Hilbert space $H$. Denote by $H_\alpha\subset H$ the subspace of all  $K^\alpha$-fixed vectors.
Denote by $P_\alpha$ the orthogonal projection to $H_\alpha$. For a double coset
$\frg\in K^\alpha\setminus \wt G/K^\beta$ we define
an operator
$$
\wt\rho(\frg):H_\beta\to H_\alpha
$$
by
$$
\wt\rho_{\alpha,\beta}(\frg):=P_\alpha \rho(g)\Bigr|_{H_\beta}
.$$

{\sc Remark.} The operator $\wt\rho(g)$ actually depends only on a double coset
containing $g$.
Indeed, for $\xi\in H_\beta$, $\eta \in H_\alpha$, and $\kappa_1\in K^\alpha$, $\kappa_2\in K^\beta$ we have
$$
\la \rho(\kappa_1 g\kappa_2)\,\xi,\,\eta \ra_{H_\alpha}=
\la \rho(g)\,\rho(\kappa_2)\,\xi,\,\rho(\kappa_1^{-1})\,\eta \ra_{H_\alpha}=
\la \rho(g)\xi,\,\eta \ra_{H_\alpha}.
$$
This expression
does not depend on $\kappa_1$, $\kappa_2$. \hfill $\boxtimes$

\sm

{\sc Remark.} Apparently, in this place we must require that the prelimit groups 
$K(\alpha)$ in (\ref{eq:GK-limits}) are compact. Otherwise
I do not see  reasons to hope for existence of non-zero fixed vectors. 
\hfill $\boxtimes$

\begin{theorem} Let $G=\GL_\fin(\infty,\Q_p)$, $K=\GL_\fin(\infty,\O_p)$.
For any $\alpha$, $\beta$, $\gamma$ and 
$$\frg_1\in 
K^\alpha\setminus \wt G/K^\beta, \qquad \frg_2\in K^\beta \setminus \wt G/K^\gamma$$
we have
\begin{equation}
\wt\rho_{\alpha,\beta}(\frg_1)\, \wt\rho_{\beta,\gamma}(\frg_2)=\wt\rho_{\alpha,\gamma}(\frg_1\circ \frg_2)
.
\label{eq:multi1}
\end{equation}
\end{theorem}

\sm

{\sc Proof.}
First, assume that the restriction of $\rho$ to $K=\GL_\fin(\infty,\O_p)$ is continuous 
in the topology of $\GL(\infty,\O_p)$. Denote by
$\ov\rho_{\alpha,\beta}(\frg)$ the following operator in $H$:
$$
\ov\rho_{\alpha,\beta}(\frg):=P_\alpha \rho(g) P_\beta,\qquad\text{where $g\in \frg$.}
$$
Representing it in a  block form
$$
H_\beta\oplus H_\beta^\bot\to H_\alpha\oplus H_\alpha^\bot
$$
we get the expression
$$
\ov\rho_{\alpha,\beta}(\frg):=\begin{pmatrix}
                         \wt\rho_{\alpha,\beta}(\frg)&0\\0&0
                        \end{pmatrix}.
$$

The statement (\ref{eq:multi1}) is equivalent to
\begin{equation}
\ov\rho_{\alpha,\beta}(\frg_1) \,\ov\rho_{\beta,\gamma}(\frg_2)=\ov\rho_{\alpha,\gamma}(\frg_1\circ \frg_2).
\label{eq:multi2}
\end{equation}
We have
\begin{multline*}
\ov\rho_{\alpha,\beta}(\frg_1)\, \ov\rho_{\beta,\gamma}(\frg_2)
=
P_\alpha\,\rho(g_1)\, P_\beta\,\rho(g_2)\,P_\gamma=
P_\alpha\,\rho(g_1)\, \Bigl(\lim_{j\to\infty} \rho(\theta_j^{[\beta]})\Bigr)\,\rho(g_2)\,P_\gamma=\\=
\lim_{j\to\infty}  P_\alpha\,\rho(g_1)\, \rho(\theta_j^{[\beta]})\,\rho(g_2)\,P_\gamma=
\lim_{j\to\infty}  P_\alpha\, \rho(g_1\theta_j^{[\beta]} g_2)\,P_\gamma
\end{multline*}
(here $\lim_{j\to\infty}$ denotes a weak limit). The sequence $g_1\theta_j^{[\beta]} g_2$
is eventually constant and we get the desired expression
$$
P_\alpha\, \rho(g_1 \circledcirc g_2)\,P_\gamma=\ov\rho(\frg_1\circ\frg_2).
$$

Next, let $\rho$ be arbitrary. The group $\GL(m,\O_p)$ centralizes $\GL^{[m]}(\infty,\O_p)$, therefore
$H_m$ is $\GL(m,\O_p)$-invariant. For $n>m$ the space $H_n$  is invariant with respect
to $\GL(n,\O_p)$ and therefore it is invariant with respect to a smaller subgroup $\GL(m,\O_p)$.
Hence $\cup_{j=0}^\infty H_j$ is invariant with respect to $\GL(m,\O_p)$. 
This is valid for all $m$, so the subspace is invariant with respect to the inductive limit $\GL_\fin(\infty,\O_p)$.
Thus we get a unitary representation of $\GL_\fin(\infty,\O_p)$ in the closure $H_*$ of
$\cup_{j=0}^\infty H_j$.  By Lemma \ref{l:adm}, this representation is continuous in the topology
of $\GL(\infty,\O_p)$, and we come to the previous case.

In $H_*^\bot$ we have no $\GL(m,\O_p)$-fixed vectors, and the statement is trivial.
\hfill $\square$.

\sm

Crucial point here is Lemma \ref{l:1}.
This picture is parallel to real classical groups and symmetric groups \cite{Olsh-semigr}, \cite{Olsh-symm},
\cite{Olsh-Howe}, \cite{Ner-book},
\cite{Ner-symmetric}, \cite{Ner-classic}, \cite{Ner-char}.
 A further discussion of $p$-adic case  is
contained in \cite{Ner-p}. 

\sm

{\sc Remark.}
It can be shown that 
in the $p$-adic case functions $\wt\rho_{\alpha,\beta} (\frg)$  do not separate
elements of $K^\alpha\setminus G/K^\beta$. Similar phenomenon is known
 for finite fields, see \cite{Olsh-semigr}.
 \hfill $\square$

\sm

{\sc Remark.}
Lemma \ref{l:1} was formulated in \cite{Ner-p} as Corollary 6.4 but its proof is incomplete
due to an incorrect definition of topology on $\GL(\infty,\O_p)$.

\sm

\section{Proof of Lemma \ref{l:1}}

\COUNTERS

{\bf\punct The symmetric group.} Denote by $S(\infty)$ the group of all permutations of the set $\N$ of natural numbers.
It has a structure of a totally disconnected  topological group defined by the condition: stabilizers of finite subsets
in $\N$ form a neighborhood basis of open subgroups in $S(\infty)$. 
Denote by $S^{[m]}(\infty)$ the group stabilizing points $1$, \dots, $m$. Clearly, open
subgroups $S^{[m]}(\infty)$ form a basis of neighborhoods of the unit in $S(\infty)$.

\sm

{\sc Remark.} This is a unique separable topology on $S(\infty)$ compatible with the group
structure.  Recall that a {\it Polish group} is a topological group, which is homeomorphic
to a complete separable metric space. There is a collection of statements about
rigidity of a choice of a Polish topology on a group, see e.g. \cite{Gao}, Section 3.2.
For instance, if two Polish topologies
on a group generate the same Borel structures, then the topologies coincide.
Of course, additive groups of all separable Banach spaces (they are Polish groups) are isomorphic as abstract groups.
But existence of such isomorphisms requires an application of the choice axiom and isomorphisms are not Borel.
\hfill $\boxtimes$

\sm

For a countable set $\Omega$ we denote by $S(\Omega)$ the group of all permutations of $\Omega$, of course
$S(\Omega)\simeq S(\infty)$.

\sm

{\bf\punct Induced representations.} Let $G$ be a totally disconnected group acting transitively on a countable set $X$,
let $R$ be a stabilizer of a point $x_0$, $\nu$ a unitary representation of $R$
in the Hilbert space $H$. Then we can define
{\it induced representation} $\mathrm{Ind}_R^G(\nu)$ of the group $G$  in the usual way 
(see, e.g., \cite{Kir}, \S 13).
Namely, we consider the space $G\times H$. Denote by $B$ its quotient with respect to the equivalence relation
$$
(x,r)\simeq (xr, \rho(r^{-1}) h),\qquad \text{where $r$ ranges in $R$.}
$$ 
Then we have a 'fiber  bundle' $B\to X=G/R$ whose fibers $H_x$ are copies of the space $H$. Transformations
$(x,h)\mapsto (gx,h)$ induce transformations of $B$. Now we consider the space of 'sections' $\psi$,
which send each point $x$ to a vector $\psi(x)\in H_x$. We define the inner product of a sections by
$$\la\psi_1,\psi_2\ra=\sum_x \la\psi_1(x),\psi_2(x)\ra_{H_x}.$$
In this way, we get a Hilbert space, the group $G$ acts on $B$ and therefore on the spaces of 
sections. This determines a unitary representation of $G$.

According the Lieberman theorem \cite{Lieb} (see also expositions in \cite{Olsh-semigr}, \cite{Ner-book})
any irreducible unitary  representation of $S(\infty)$ is induced from an irreducible
representation of a subgroup of the type $S(m)\times S^{[m]}(\infty)$ trivial on the factor $S(\infty-m)$.
We need  the following fact (see \cite{Ner-book}, Corollary VIII.1.5), it immediately
follows from the Lieberman theorem.

\sm

\begin{lemma}
\label{l:lieb}
 For any unitary representation $\rho$ of $S(\infty)$ the sequence $\rho(\theta_j)$
 weakly converges to the orthogonal projection to the space of vectors fixed by all operators $\rho(g)$.
\end{lemma}

\sm

{\bf \punct Oligomorphic groups.}
Recall that a closed subgroup $G$ in $S(\Omega)$ is called {\it oligomorphic}
if  it has a finite number of orbits on each finite product $\Omega\times \dots\times \Omega$.
We need the following Tsankov theorem \cite{Tsa}:

\sm

\begin{theorem}
Any unitary representation of an oligomorphic group $G$ is a {\rm(}countable or finite{\rm)} direct sum of
irreducible representations. For any irreducible representation
$\rho$ of $G$ there are open subgroups $R\subset \wt R$  such that $R$ is a normal subgroup of  finite index in $\wt R$
and
\begin{equation}
\rho=\mathrm{Ind}_{\wt R}^G(\nu),
\label{eq:Ind}
\end{equation}
where $\nu$ is an irreducible representation of $\wt R$ trivial on $R$.
\end{theorem}

\begin{corollary}
\label{cor:o}
Any irreducible representation of an oligomorphic group $G$ is a subrepresentation of a
quasiregular representation in $\ell^2$
on some homogeneous space $G/R$, where $R$ is an open subgroup in $G$.
\end{corollary}

{\sc Proof.} Let $\tau$ be a unitary representation of the group $\wt R/R$.
Denote by $\tau_{\circ}$ the same representation considered as a representation 
of $\wt R$ trivial on $R$. Denote by 
$\mathrm{Reg}$ the regular representation of $\wt R/R$. Denote by $\tau^0$ the trivial
(one-dimensional) representation of $\wt R/R$.

It is easy to see that
$$
\mathrm{Ind}_R^{\wt R}(\tau^0_\circ)=\mathrm{Reg}_\circ.
$$
Let $\nu$ be as above. Then $\nu$ is a subrepresentation of 
$\mathrm{Reg}_\circ$, therefore $\rho$ given by (\ref{eq:Ind})
is a subrepresentation of
$$
\mathrm{Ind}_{\wt R}^G\Bigl(
\mathrm{Ind}_R^{\wt R}(\tau^0_\circ)\Bigr)= \mathrm{Ind}_R^{G}(\tau^0_\circ).
$$
 The last representation is the quasiregular representation of $G$ in $\ell^2(G/R)$.
 \hfill $\square$

\sm

{\bf \punct Definitions.}

\sm

{\sc a) Modules.} Denote
by $\Z_{p^k}$
the residue rings $\Z/p^k\Z$. 
A module over $\Z_{p^k}$ is nothing but an Abelian $p$-group
whose elements have orders $\le p^k$.

The $p$-adic integers $\O_p$
are the inverse limit
\begin{equation}
\O_p=\lim_{\longleftarrow} \Z_{p_k}.
\label{eq:projectiveO}
\end{equation}
A reduction of a $p$-adic integer $x$ modulo $p^k$ we denote by
$$\l x\r_{p^k}\in \Z_{p_k}.$$
We will use the same notation for reductions of vectors and matrices.

For each $k$  define a $\Z_{p^k}$-module $V(\Z_{p^k})$ as the space of all sequences $z=(z_1,z_2,\dots)$,
where $z_j\in \Z_{p^k}$ and $z_l=0$ for sufficiently large $l$. We equip this space
with the discrete topology.

Next, we define an $\O_p$-module $V(\O_p)$ as the space of all sequences $z=(z_1,z_2,\dots)$,
where $z_j\in\O_p$ and $|z_j|\to \infty$ as $j\to\infty$. In other words,
$$
V(\O_p)= \lim_{\longleftarrow} V(\Z_{p^k}),
$$
we equip this space with the topology of a projective limit. A sequence
$z^{(l)}\in V(\O_p)$ converges if
all reductions  $\l z^{(l)}\r_{p^k} \in V(\Z_{p^k}) $ are eventually constant.
The same topology is induced by the norm
$$
\|z\|:=\max\limits_j |z_j|.
$$

We also define 'dual' modules $V^\circ(\Z_{p^k})$, $V^\circ(\O_p)$ consisting of vector-columns
satisfying same properties.

\sm

{\sc b) Groups
$\GL(\infty, \Z_{p^k})$}. We define this group as a group of all infinite matrices over $\Z_{p^k}$
such that each row and each column contains only finite number of nonzero elements.
The group $\GL(\infty, \Z_{p^k})$ acts by automorphisms on the module
$V(\Z_{p^k})\oplus V^\circ(\Z_{p^k})$ by the transformations
$$
g:\,\, (v,w^\circ)\to (vg, g^{-1} w^\circ).
$$
Thus we have an embedding to a symmetric group,
$$
\GL(\infty, \Z_{p^k})\to S\Bigl( V(\Z_{p^k})\oplus V^\circ(\Z_{p^k})\Bigr). 
$$
We equip $\GL(\infty, \Z_{p^k})$ with the induced topology.
For any collection of vectors $v_1$,\dots, $v_l\in V$ and covectors 
$w^\circ_1$, \dots, $w^\circ_n$ its stabilizer 
\begin{equation}
G(v_1,\dots,v_l;w^\circ_1, \dots, w^\circ_n  )
\label{eq:Gvw}
\end{equation}
is an open subgroup in $\GL(\infty, \Z_{p^k})$. By definition,
such subgroups form a basis of neighborhoods of unit in our group.

\sm

Next, for each $m$ we define the subgroup $\GL^{[m]}(\infty,\Z_{p^k})\subset \GL(\infty,\Z_{p^k})$
consisting of all matrices of the form 
$
\begin{pmatrix}
 1_m&0\\
 0&*
\end{pmatrix}
$. This group has the form $G(e_1,\dots,e_m; f^\circ_1,\dots, f^\circ_m)$, where $e_j$ is the standard basis 
in $V$ and $f_j^\circ$ is the standard basis in $V^\circ$. 
Since actually vectors and covectors $v_i$ and $w_j$ in (\ref{eq:Gvw}) have only finite
number of nonzero coordinates, each stabilizer $G(\dots)$ contains some group $\GL^{[m]}(\infty,\Z_{p^k})$.

Thus, {\it the subgroups $\GL^{[m]}(\infty,\Z_{p^k})$  form  a basis of neighborhoods of unit in our group.}

We can also define the topology in the following way. A sequence $g_l$ converges to $g$
if for each $i$ the sequence of $i$-th rows (resp. columns) of $g_l$ coincides with the $i$-th row
(resp. column) of $g$
for sufficiently large $l$.

\sm

{\sc c) The group $\GL(\infty,\O_p)$.} We have natural homomorphisms of rings 
$\Z_{p^k}\to \Z_{p^{k-1}}$ and therefore homomorphisms of groups
$$
\GL(\infty,\Z_{p^k})\to \GL(\infty,\Z_{p^{k-1}}).
$$
We define the group $\GL(\infty,\O_p)$ as the projective limit
$$
\GL(\infty,\O_p):=\lim_{\longleftarrow} \GL(\infty,\Z_{p^k}).
$$
In other words, this group consists of all infinite matrices $g$ over $\O_p$
such that $\l g\r_{p^k}\in \GL(\infty,\Z_{p^k})$ for all $k$.

\sm

The topology on $\GL(\infty,\O_p)$ is the topology of the projective limit.
A sequence $g^{(j)}$ converges to $g$ if
$\l g^{(j)}\r_{p^k}\in \GL(\infty,\Z_{p^k})$ converges to  $\l g\r_{p^k}$
for all $k$.

\sm

{\sc d) Open subgroups in $\GL(\infty,\O_p)$.}
For nonnegative integers $m$, $k$ we define subgroups
$\GL_k^{[m]}(\infty,\O_p)$ consisting of $(m+\infty)$-block matrices having the form
$$
\begin{pmatrix}
1+p^k A& p^k B\\
p^k C& D
\end{pmatrix},
$$
where $A$, $B$, $C$, $D$ are matrices over $\O_p$. These subgroups are open
 and form a basis of neighborhoods of 1. 

\sm

We define a {\it congruence subgroup} $\GL_k(\infty,\O_p)$ in $\GL(\infty,\O_p)$
as the  subgroup consisting of matrices
having the form $1+p^k Q$, where $Q$ is a matrix over $\O_p$ (congruence subgroups are not open).

\sm

{\bf\punct  Lemmas.} 
Next, we apply the following general statement, see \cite{Ner-book}, Proposition VII.1.3.

\sm

\begin{proposition}
\label{pr:}
Let $G$ be a topological group, $G_1\supset G_2\supset\dots$ be a sequence of subgroups 
such that any neighborhood of unit in $G$ contains a subgroup $G_j$. Let $\rho$ be 
a unitary representation of $G$ in a Hilbert space $H$. Denote by 
$H_k$ the space of vectors invariant with respect to $G_k$. Then $\cup H_k$
is dense in $H$.
\end{proposition}

\begin{corollary}
 Any unitary representation $\rho_j$ of $\GL(\infty,\O_p)$ can be decomposed
 as a direct sum $\oplus_{k=1}^\infty \rho_k$, where
 $\rho_k$ is trivial on the congruence subgroup $\GL_k(\infty,\O_p)$.
\end{corollary}

{\sc Proof.}
We apply Proposition \ref{pr:} to the group $\GL(\infty,\O_p)$ and  
the sequence of congruence subgroups $\GL_k(\infty,\O_p)$.
Since a subgroup $\GL_k(\infty,\O_p)$ is normal, for $h\in H_k$,
$g\in \GL(\infty,\O_p)$, $r\in \GL_k(\infty,\O_p)$, we have
$$
\rho(r)\, \rho(g)\, h= \rho(g)\, \rho(g^{-1} r g)\, h.
$$
Since the congruence subgroup is normal, $g^{-1} r g\in \GL_k(\infty,\O_p)$,
therefore  
$$\rho(g^{-1} r g) h=h,$$ i.e.,  $h\in H_k$.
Therefore the subspace $H_k$ is invariant with respect
to the whole group $\GL(\infty,\O_p)$ and the congruence subgroup
acts in $H_k$ trivially.
Thus, 
$$
H=\oplus_{k=1}^\infty (H_{k}\ominus H_{k-1}).
$$
In each $H_{k}\ominus H_{k-1}$ we have an action of $\GL(\infty,\Z_{p^k})$.

Thus, it is sufficient to prove Lemma \ref{l:1} for groups $\GL(\infty,\Z_{p^k})$.

\sm

Recall that we can consider the group $S(\infty)$
as a group of 0-1-matrices.

\begin{lemma}
\label{l:gen}
 For any $m$ the group $\GL(\infty,\Z_{p^k})$ is generated by the
 subgroups
 $S(\infty)$ and $\GL^{[m]}(\infty,\Z_{p^k})$.
\end{lemma}

{\sc Proof.} Consider the subgroup $G$ generated by these subgroups.
Clearly, $G$ contains all
groups $\GL(n,\O_p)$.  Indeed, for $y\in \GL(n,\O_p)$  for 
$N>\max (m,n)$ we have $\theta_N y\theta_N^{-1}\in \GL^{[m]}(\infty,\Z_{p^k})$.

Fix $g\in \GL(\infty,\Z_{p^k})$.
For suffitiently large $\beta$ 
the expression of $g$ as a block $(m+\beta+\infty)$-matrix has the form
$$
g=\begin{pmatrix}
   g_{11}&g_{12}&0\\
   g_{21}&g_{22}&g_{23}\\
   0&g_{32}&g_{33}
  \end{pmatrix}.
$$
Multiplying this matrix from the right by an appropriate matrix of the form
\begin{equation}
\begin{pmatrix}
r_{11}&r_{12}&0\\
r_{12}&r_{22}&0\\
0&0&1
\end{pmatrix}
\label{eq:r}
\end{equation}
we can obtain a matrix of the form
\begin{equation}
g'=\begin{pmatrix}
   1&0&0\\
   g'_{21}&g'_{22}&g'_{23}\\
   0&g'_{32}&g'_{33}
  \end{pmatrix}.
  \label{eq:g-prime}
\end{equation}
Indeed, we can regard rows $u_1$, \dots, $u_m$ of the matrix
$
\begin{pmatrix}
   g_{11}&g_{12}
\end{pmatrix}
$
as elements of the module $\Z_{p^k}^{m+\beta}$. Since the matrix
$g$ is invertible, 
  the  matrix $\l g\r_{p}$ over the finite field
 $\Z_p$ is invertible.
Therefore  the matrix
$\l\begin{pmatrix}
   g_{11}&g_{12}
\end{pmatrix}\r_p$
 is nondegenerate. This implies that rows of
$u_1$, \dots, $u_m$ generate a submodule isomorphic $\Z_{p^k}^m$.
Adding an appropriate collection $v_1$, \dots, $v_\beta$ 
we can obtain a basis of the module $\Z_{p^k}^{m+\beta}$.
The matrices (\ref{eq:r}) determines automorphisms of $\Z_{p^k}^{m+\beta}$.
We send $u_1$, \dots, $u_m$, $v_1$, \dots, $v_\beta$
to the standard basis in $\Z_{p^k}^{m+\beta}$.

Thus we came to a matrix $g'$ given by (\ref{eq:g-prime}).
Multiplying  $g'$ from the left by
by
$$
\begin{pmatrix}
 1&0&0\\
 -g_{21}'&1&0\\
 0&0&1
\end{pmatrix}
$$
we come to
$$
 \qquad\qquad\qquad\qquad 
g''=\begin{pmatrix}
   1&0&0\\
   0&g'_{22}&g'_{23}\\
   0&g'_{32}&g'_{33}
\end{pmatrix}\in \GL^{[m]}(\infty,\Z_{p^k}).
 \qquad\qquad\qquad\qquad \square
$$

\begin{lemma} 
Groups $\GL(\infty, \Z_{p^k})$ are oligomorphic.
\end{lemma}

{\sc Proof.}
We have an action of $\GL_\infty(\infty,\Z_{p^k})$ on $V(\Z_{p^k})\oplus V^\circ(\Z_{p^k})^n$,
i.e. on collections $(v_1,\dots, v_n;\,w_1^\circ,\dots,w^\circ_n)$  of vectors and covectors.
We must show that there is an finite set containing representatives of all $\GL_\infty(\infty,\Z_{p^k})$-orbits.
Denote by $V_N$ the submodule in $V(\Z_{p^k})$ consisting of vectors whose coordinates with numbers $>N$
vanish, $V_N=\Z_{p^k}^N$. A block $(N+\infty)$-matrix of the form $g=\begin{pmatrix} a&0\\0&1\end{pmatrix}$ 
induces an automorphism of $V_N$. We can sent $v_1$, \dots, $v_n$ to the submodule $V_n\subset V_N$.

Next consider the action of the group $\GL^{[n]}(\infty,\Z_{p^k})$ on collections of vectors and covectors.
It does not change  vectors and  first $n$ coordinates of covectors. The same argument as above shows that we
can put all covectors to the module $V^\circ_{2n}$.

Thus any orbit intersects a finite set $V(\Z_{p^k})^n\oplus V^\circ (\Z_{p^{k}})^{2n}$
\hfill $\square$

\sm

{\bf \punct Proof of Lemma  \ref{l:1}.} By Corollary \ref{cor:o}, it is sufficient to prove 
the statement for a quasiregular representations $\rho$ of $\GL(\infty,\Z_{p^k})$
in a space $\ell^2(X)$, where $X=\GL(\infty,\Z_{p^k})/P$ is coset space with respect
to an open subgroup $P$.
For $x\in X$
denote by $\delta_x$ an element of $\ell^2$, which is 1 at $x$ and 0 at other points.

By Lemma \ref{l:lieb}, the weak limit of $\rho(\theta_j)$ exists and coincides with orthogonal projection
to $S(\infty)$-fixed vectors. Let 
$$\psi=\sum_{x\in X} c_x \delta_x\ne 0$$
be such a vector. 
For $\sigma\in S(\infty)$ we have $c_{\sigma x}=c_x$. If $x$ is not fixed by $S(\infty)$,
then its orbit is infinite. Since $\psi\in \ell^2$, we have $c_x=0$.
Thus, $\psi$ has the form
$$
\psi=\sum_{\text{$x$: $\sigma x=x$ for all $\sigma\in S(\infty)$ }} c_x\delta_x.
$$
A stabilizer  of $x$ is an open subgroup in $\GL(\infty,\O_p)$.
It contains some subgroup of the form $\GL^{[m]}(\infty,\O_p)$.
On the other hand it contains $S(\infty)$. By Lemma
\ref{l:gen}, the stabilizer of $x$ is the whole group $\GL(\infty,\O_p)$.
Thus, the space $X$ consists of one point.
This finishes a proof.
\hfill $\square$

\section{Admissibility}

\COUNTERS

Here we prove Lemma \ref{l:adm}.
It is sufficient to prove $\Leftarrow$, since
the implication $\Rightarrow$ immediately follows from Proposition \ref{pr:}.

\sm

{\bf\punct A normal form for double cosets.%
\label{ss:normal}}

\begin{lemma}
\label{l:long}
{\rm a)} Any double coset of $\GL_\fin(\infty,\O_p)$ with respect to $\GL^{[m]}_\fin(\infty,\O_p)$
 contains an element of $\GL(3m,\O_p)$.
 
 \sm
 
{\rm b)} The same statement holds for double cosets of $\GL(\infty,\O_p)$ by $\GL^{[m]}(\infty,\O_p)$. 

\sm

{\rm c)} The natural map
\begin{multline}
\GL^{[m]}_\fin(\infty,\O_p)\setminus \GL_\fin(\infty,\O_p)/\GL^{[m]}_\fin(\infty,\O_p)
\longrightarrow\\\longrightarrow
\GL^{[m]}(\infty,\O_p)\setminus \GL(\infty,\O_p)/\GL^{[m]}(\infty,\O_p)
\label{eq:map-cosets}
\end{multline}
is a bijection.

\sm

{\rm d)}
Let $M\ge 3m$. Let for $g_1$, $g_2\in \GL(M,\O_p)$
there are $q$, $r\in \GL^{[m]}(\infty,\O_p)$
such that $g_1=q g_2 r$. Then for sufficiently large $N$
depending only on $M$ there are 
$$\text{$q'$, $r'\in \GL(N,\O_p)\cap \GL^{[m]}(\infty,\O_p)$}$$
such that $g_1=q' g_2 r'$.
%
\end{lemma}

{\sc Proof. a), b)} In the both cases we can apply the following reduction.
Represent $g\in\GL_\fin(\infty,\O_p)$ as a block matrix 
$
g=\begin{pmatrix}
   g_{11}&g_{12}\\
   g_{21}&g_{22}
  \end{pmatrix}
$
of size $(m+\infty)$.
Multiplying it from the right by matrices of the form $\begin{pmatrix}
                                                        1&0\\
                                                        0&u
                                                       \end{pmatrix}\in \GL^{[m]}$
we can reduce it to the $(m+m+\infty)$-block form
$$
g'=\begin{pmatrix}
   g'_{11}&g'_{12}&0\\
   g'_{21}&g'_{22}&g'_{23}\\
   g'_{31}&g'_{32}&g'_{33}
  \end{pmatrix}
$$
(in fact $g'_{12}$ can be made lower-triangular).
Applying similar left multiplication we can make $g'_{31}=0$.

\sm

Next, we multiply $g'$ from the left and from the right by elements
of $\GL^{[2m]}_\fin$ to simplify $g'_{33}$ (such multiplications do
not change blocks $g'_{11}$, $g'_{12}$, $g'_{21}$, $g'_{22}$).
If the reduction $\l g'_{33}\r_p$  is nondegenerate, we can
make $\l g'_{33}\r_p=1$ and $g'_{33}=1$. 

However $\l g'_{33}\r_p$ can be degenerate
and 
$$\dim\ker g'_{33}=\dim\coker g'_{33}:=\gamma\le m.$$ 
In this case, we can 
make from $\l g'_{33}\r_p$ a matrix $\begin{pmatrix} 0_\gamma&0\\0&1_\infty
                                  \end{pmatrix}$
                                    and reduce $g'_{33}$ to the form
 $\begin{pmatrix} pA&pB\\pC&1_\infty+pD
                                  \end{pmatrix}$,
                                  where $A$, $B$, $C$, $D$ are  matrices over $\O_p$.
Applying a right multiplication by
$
\begin{pmatrix}
 1_{2m+\gamma}&0\\0&1+pD
\end{pmatrix}^{-1}
\in
\GL^{[2m+\gamma]}\subset \GL^{[2m]}
$,
we 'kill' $pD$ and come to 
to a block $(m+m+\gamma+\infty)$-matrix 
 of the form
 $$
  g'':=
 \begin{pmatrix}
g_{11}''&g_{12}''&0&0\\
g_{21}''&g_{22}''&g_{23}''&g_{24}''\\
0&g_{32}''&g_{33}''&g_{34}'\\
0&g_{42}''&g_{43}'&1\\  
 \end{pmatrix}
 .
 $$
 Multiplying by 
 the matrix
 $$
 \begin{pmatrix}
  1&0&0&0\\
  0&1&0&-g_{24}''\\
  0&0&1&-g_{34}''\\
  0&0&0&1
 \end{pmatrix}
 $$
 from the left, we kill $g_{24}''$, $g_{34}''$ (and change only $g_{22}''$, $g_{23}''$, $g_{33}''$, $g_{34}''$).
 In the same way (by a multiplication from the right) we kill $g_{42}''$, $g_{43}''$.

\sm

{\sc d)}
Denote $l:=M-m$. We wish to verify the following statement:
if for given $g_1$, $g_2\in \GL(m+l,\O_p)$  there exist $\xi$, $\eta\in \GL^{[m]}(\infty,\O_p)$
satisfying the equation
\begin{equation}
g_1 \xi=\eta g_2,
\label{eq:gg}
\end{equation}
then there exist $\xi'$, $\eta'\in \GL^{[m]}_\fin(\infty,\O_p)$ satisfying the same equation.
Let us write the equation (\ref{eq:gg}) as a condition for $(m+l+\infty)$-block matrices,
\begin{equation}
\begin{pmatrix}
 a&b&0\\
 c&d&0\\
 0&0&1_\infty
\end{pmatrix}
\begin{pmatrix}
 1_m&0&0\\
 0&x&y\\
 0&z&u
\end{pmatrix}
=
\begin{pmatrix}
 1_m&0&0\\
 0&X&Y\\
 0&Z&U
\end{pmatrix}
\!\!
\begin{pmatrix}
 A&B&0\\
 C&D&0\\
 0&0&1_\infty
\end{pmatrix}
\label{eq:}
\end{equation}
(the matrices in the left hand side denote $g_1$ and $\xi$, the matrices in the right hand side $\eta$ and $g_2$)
or
\begin{equation}
\begin{pmatrix}
 a&bx&by\\
 c&dx&dy\\
 0&z&u
\end{pmatrix}=
\begin{pmatrix}
 A&B&0\\
 XC&XD&Y\\
 ZC&ZD&U
\end{pmatrix}.
\label{eq:system}
\end{equation}
Let $\varkappa$ be an infinite invertible matrix over $\O_p$. Then
the transformations
$$
 Z\mapsto \varkappa z,\qquad U\mapsto\varkappa U, \qquad z\mapsto \varkappa z, \qquad u\mapsto\varkappa u
$$
send a solution of the system of equations (\ref{eq:system}) to a solution.
We can find  a new solution, where
$Z$ has a form $\begin{pmatrix} Z'\\0\end{pmatrix}$,
the size of this matrix is $(l+\infty)\times l$.
By (\ref{eq:system}) the new $z$ is $ZD=\begin{pmatrix} Z'D\\0\end{pmatrix}$.
Applying a similar transformation
$$
y\mapsto y\lambda,\qquad u\mapsto u\lambda,\qquad  Y\mapsto Y\lambda,\qquad  U\mapsto u\lambda
$$
we can get a solution of system (\ref{eq:system}) with $y$ and $Y$ having form $\begin{pmatrix}*&0\end{pmatrix}$.
Thus we have a solution of the equation (\ref{eq:})
with finitary indeterminates $z$, $Z$, $y$, $Y$.
Now the indeterminant factor in the left hand side of (\ref{eq:}) can be written in the 
$(m+l+l+\infty)$-block form
\begin{equation}
\begin{pmatrix}
 1_m&0&0&0\\
 0&x'&y'&0\\
 0&z'&u'_{11}&u'_{12}\\
 0&0&u'_{21}&u'_{22}
\end{pmatrix}.
\label{eq:4}
\end{equation}
The only equation in (\ref{eq:system}) containing $u$ is $u=U$.
Since the matrices $\l\begin{pmatrix}
                       x'&y'
                      \end{pmatrix}\r_p$
                      and 
                      $\l\begin{pmatrix}
                       x'\\z'
                      \end{pmatrix}\r_p$
                      are nondegenerate,
                      we can choose  $u_{11}'$ such that a matrix 
                      $\l\begin{pmatrix}
                       x'&y'\\z'&u_{11}'
                      \end{pmatrix}\r_p$ also is nondegenerate.
 We set $u_{12}'=0$, $u_{21}'=0$, $u_{22}'=1_\infty$. Then the matrix (\ref{eq:4}) is invertible.                     
 
 Finally, we found new $U$ from the equation $u=U$.
Then 3 factors in (\ref{eq:}) are invertible and therefore the fourth factor also is invertible.
A finitary solution of (\ref{eq:}) is obtained. Actually,
$$
\xi',\,\eta' \in \GL^{[m]}(\infty,\O_p)\cap \GL(m+2l,\O_p).
$$

\sm

{\sc c)} The surjectivity follows from a) and b). 
The injectivity follows from d).
\hfill $\square$

\sm

{\bf\punct The metric on the space of double cosets.%
\label{l:metric}} Here we prove the following lemma:

\begin{lemma}
The maps {\rm (\ref{eq:map-cosets})} are homeomorphisms.
\label{l:continuity}
\end{lemma}

Fix $m$.
For each $M\ge 3m$ we have a natural partition of the group
$\GL(M,\O_p)$ into subsets $\GL(M,\O_p)\cap \frg$, where 
$\frg$ are double cosets of $\GL(\infty,\O_p)$
with respect to $\GL^{[m]}(\infty,\O_p)$. Denote by $\cK_M$ the quotient space,
 According 
Lemma \ref{l:long}.d, elements of partitions are compact, therefore quotients are compact.
For $M<M'$ the natural map $\cK_M\to \cK_{M'}$ is continuous,
By Lemma  \ref{l:long}.a, it is a bijection, hence it is a homeomorphism.
This also implies that the bijections
$$
\cK_M\longleftrightarrow \GL^{[m]}_\fin(\infty,\O_p)\setminus \GL_\fin(\infty,\O_p)/\GL^{[m]}_\fin(\infty,\O_p)
$$
are homeomorphisms.
Also it is clear that the maps 
$$
\cK_M\longrightarrow \GL^{[m]}(\infty,\O_p)\setminus \GL(\infty,\O_p)/\GL^{[m]}(\infty,\O_p)
$$
are continuous. We must proove a continuity of the inverse map.

\sm

Define a left-right-invariant metric on $\GL(\infty,\O_p)$ by
$$
d(z,u)=\max_{i,j}|z_{ij}-u_{ij}|.
$$

{\sc Remark.} This metric determines on each group $\GL(n,\O_p)$
the standard topology.
On the whole group $\GL(\infty,\O_p)$ it determines a  nonseparable topology,
which 
is stronger than the natural topology.
Restriction of the metric to $\GL_\fin(\infty,\O_p)$ induces a topology that is weaker than the natural topology.
\hfill $\square$

\sm

Recall that the Hausdorff metric on the space of compact subsets of a metric space
is given by the formula
$$
\dist_H (A,B):=\max\Bigl[ \max_{x\in A} \min_{y\in b} d(x,y), \, \max_{y\in b} \min_{x\in A} d(x,y) \Bigr].
$$
Restricting this metric to elements of the partition of $\GL(M,\O_p)$
we get a metric on the double coset space
\begin{equation}
\label{eq:distH}
\dist_H^M(\frg_1, \frg_2)=\dist_H \Bigl(\frg_1\cap \GL(M,\O_p), \frg_2\cap \GL(M,\O_p)\Bigr)
\end{equation}
 compatible with the topology on $\cK_M$.

Next, define another metric on $\GL^{[m]}(\infty,\O_p)/\GL(\infty,\O_p)/\GL^{[m]}(\infty,\O_p)$.
Let $\frg_1$, $\frg_2$ be double cosets. Fix
$g\in \frg_1$. Then 
\begin{equation}
\dist (\frg_1,\frg_2)=\inf_{z\in \frg_2} d(g,z)
\label{eq:dist}
\end{equation}
(the result does not depend on $g$).

\begin{lemma}
 These metrics coincide.
\end{lemma}

We have an obvious inequality
$$
\dist_H^{3m} (\frg_1,\frg_2)\ge 
\dist (\frg_1,\frg_2).
$$
The 
inverse inequality follows from the 
following lemma:

 \begin{lemma}
 \label{l:dist}
 Let $\frg_1$, $\frg_2$ be  double cosets.
  Let $g_1\in \frg_1$, $g_2\in \frg_2$. Let $u\in \frg_1\cap \GL(3m,\O_p)$.
Then there exist $w\in \frg_2\cap \GL(3m,\O_p)$ such that 
$$
d(u,w)\le d(g_1,g_2).
$$
 \end{lemma}

{\sc Proof .} Let $d(g_1,g_2)=p^{-k}$. Let $u= q g_1 r$, where $q$, $r\in \GL^{[m]}(\infty,\O_p)$.
We take $v:=q g_2 r$. Then we can make a reduction as in the proof of Lemma \ref{l:long}
using only elements of the congruence subgroup, i.e. we found 
$w=t v s\in \GL(3m,\O_p)$ with $t$, $s\in \GL^{[m]}(\infty,\O_p)\cap \GL_k(\infty,\O_p)$.
Then we have $d(v,u)=p^{-k}$, $d(u,w)\le p^{-k}$.
\hfill $\square$

\sm

{\sc Proof of Lemma \ref{l:continuity}} It sufficient to show that for $\frg$ for each $k$ there is 
a neighborhood $\cN$ of $\frg$ in the sense of 
$\GL^{[m]}(\infty,\O_p)\setminus \GL(\infty,\O_p)/\GL^{[m]}(\infty,\O_p)$
such that for each $\frh\in \cN$ we have $\dist(\frg,\frh)\le p^{-k}$.

Choose $g\in\frg$, choose $N$ such that  $\l g\r_{p^k}$
 has the following $(m+N+\infty)$-block form
 $$\l g\r_{p^k}=\begin{pmatrix}
     u_{11}&u_{12}&0\\u_{21}&u_{22}&u_{23}\\
     0&u_{32}&u_{33}
    \end{pmatrix}.
    $$
    Next, we consider an open subgroup $\GL^{[m+N]}_k(\infty,\O_p)$ and consider
    the neighborhood $\cO:=g \,\GL^{[m+N]}_k(\infty,\O_p)$ of $g$. Let $r\in\GL^{[m+N]}_k(\infty,\O_p)$,
    then
    $h=g r\in \cO$. 
    Consider the matrix $\l r\r_{p^k}\in \GL^{[m+N]}(\infty,\Z_{p^k})$. 
    Let us regard it  as a matrix $\wt r\in \GL^{[m+N]}(\infty,\O_p)$ composed
    of $p$-adic integers  contained in the set $0$, $1$, \dots, $p^{k-1}$.
    Consider the matrix
    $gr\wt r^{-1}$, which is contained in the same double coset. 
    Then $\l gr\wt r^{-1}\r_{p^k}=\l g\r_{p^k}$.
    Thus,
    $$
    |g-gr\wt r^{-1}|\le p^{-k}.
    $$
  We apply Lemma \ref{l:dist} and this completes the proof. \hfill $\square$
  
  \sm

{\bf \punct End of proof of Lemma \ref{l:adm}.%
\label{ss:matrix-elements}}
It is sufficient to show that matrix elements of the form
$$
\la \rho(g) \xi,\eta\ra_H, \qquad\text{where $\xi$, $\eta\in \cup H_j$}
$$
have continuous extensions to the whole group $\GL(\infty,\O_p)$. We can assume that 
$\xi$, $\eta\in  H_m$. Such matrix elements  are continuous functions 
 on the inductive limit $\GL_\fin(\infty,\O_p)$, which are constant on double cosets with respect
 to $\GL^{[m]}_\fin(\infty,\O_p)$.
By Lemma \ref{l:continuity} they are continuous on $\GL(\infty,\O_p)$.

 \tt
 
 \noindent
 Math. Dept., University of Vienna;c/o Wolfgang Pauli Institute \\
 Institute for Theoretical and Experimental Physics (Moscow); \\
 MechMath Dept., Moscow State University;\\
 Institute for Information Transmission Problems.\\
 URL: http://mat.univie.ac.at/$\sim$neretin/

\end{document}